\documentclass[12pt, twoside]{article}
\usepackage{somedefs,tracefnt,latexsym,,rawfonts,epsfig,amsfonts,amssymb,latexsy
m,enumerate,amscd}

\def\Bbb{\mathbb}
\def\n{\noindent}
\def\ms{\medskip}
\def\bi{\bibitem}
\def\nb{\newblock}

\title{Topology of $3$-manifolds and a class of groups
II\footnote{accepted for publication  
in the special issue of the {\it Boletin de la Sociedad Matematica 
Mexicana} dedicated to Professor Francisco Gonzalez Acuna on the occasion
of his 60'th birthday.}}
\author{S. K. Roushon}
\date{March 01, 2004}

\begin{document}
\pagestyle{myheadings}
\markboth{\centerline{Topology of $3$-manifolds and a class of groups
II}$\!\!\!\!\!\!$}{\centerline{S.K. Roushon}$\!\!\!\!\!\!$}
\maketitle

\begin{abstract}
This is a continuation of an earlier preprint (\cite{R2}) under the same
title. These papers grew out of an attempt to find a suitable finite
sheeted covering of an aspherical $3$-manifold so that the cover either
has infinite or trivial first homology group. With this motivation we
defined a new class of groups. These groups are in some sense eventually
perfect. Here we prove results giving several classes of examples of
groups which do (not) belong to this class. Also we prove some basic
results on these groups and state two conjectures. A direct application of
one of the conjectures to the virtual Betti number conjecture 
is mentioned. For completeness, here we reproduce parts of \cite{R2}.   

\ms
\n

\n
\end{abstract}
\baselineskip 14pt
\newpage
\setcounter{section}{-1}
\section{Introduction}

\newenvironment{proof}{{\it Proof.}}{$\hfill{\Box}$}

The main motivation to this paper and \cite{R2} came from $3$-manifold
topology while trying to find a suitable finite sheeted covering of an
aspherical $3$-manifold so that the cover has either infinite or trivial
first integral homology group. In \cite{R} it was proved that $M^3\times
{\Bbb D}^n$ is topologically rigid for $n>1$ whenever $H_1(M^3, {\Bbb Z})$ 
is infinite. Also the same result is true  when $H_1(M^3, {\Bbb Z})$ is
$0$. The remaining case is when $H_1(M^3, {\Bbb Z})$ is nontrivial
finite. There are induction techniques in surgery theory which can be
used to prove topological rigidity of a manifold if certain finite
sheeted coverings of the manifold are also topologically rigid.
In the case of manifolds with nontrivial finite first integral homology
groups there is a natural finite sheeted cover, namely, the one which
corresponds to the commutator subgroup of the fundamental group.

So we start with a closed aspherical $3$-manifold $M$ with nontrivial
finite first integral homology group and consider the finite sheeted
covering $M_1$ of $M$ corresponding to the commutator subgroup. If
$H^1(M_1, {\Bbb Z})\neq 0$ or $H_1(M_1, {\Bbb Z})=0$ then we are done 
otherwise we again take the finite sheeted cover of $M_1$ corresponding 
to the commutator subgroup and continue. The group theoretic
conjecture (Conjecture 0.2) in this article implies that this process
stops in the sense that for some $i$ either $H^1(M_i, {\Bbb Z})\neq 0$ or 
$H_1(M_i, {\Bbb Z})=0$. 

Motivated by the above situation we define the following class of groups. 

\ms
\n
{\bf Definition (0.0).} An abstract group $G$ is called {\it adorable} if
$G^i/G^{i+1} = 1$ for some $i$, where $G^i=[G^{i-1}, G^{i-1}]$, the
commutator subgroup of $G^{i-1}$, and $G^0=G$.  The smallest $i$ for which
the above property is satisfied is called the {\it degree of adorability}
of $G$. We denote it by $doa(G)$.

\ms
Obvious examples of adorable groups are finite groups, perfect groups,
simple groups and solvable groups. The second and third class of groups 
are adorable groups of degree $0$. The free products of perfect groups are
adorable (in fact perfect). The nontrivial abelian groups and symmetric
groups on $n\geq 5$ letters are adorable of degree $1$. Another class of
adorable groups are $GL(R)=Lim_{n\to \infty}GL_n(R)$. Here $R$ is any ring
with unity and $GL_n(R)$ is the multiplicative group of $n\times n$
invertible matrices. These are adorable groups of degree $1$. This follows
from the Whitehead lemma which says that the commutator subgroup of
$GL_n(R)$ is generated by the elementary matrices and the group generated
by the elementary matrices is a perfect group. Also $SL_n({\Bbb C})$, the
multiplicative group of $n\times n$ matrices with complex entries is a
perfect group. In fact we will prove that any connected Lie group is
adorable as an abstract group. The full braid groups on more than $4$
strings are adorable of degree $1$. 

We observe the following two elementary facts in the next section.

\ms
\n
{\bf Theorem (1.8).} {\it A group $G$ is adorable if and only if
there is a filtration $G_n< G_{n-1}<\cdots <G_1<G_0=G$ of
$G$ so that $G_i$ is normal in $G_{i-1}$, $G_{i-1}/G_i$ is abelian for
each $i$ and $G_n$ is a perfect group.}

\ms
\n
{\bf Theorem (1.10).} {\it Let $H$ be a normal subgroup of an adorable
group $G$. Then $H$ is adorable if one of the following conditions is 
satisfied.
\begin{itemize}
\item $G/H$ is solvable.
\item for some $i$, $G^i/H^i$ is abelian.
\item for some $i$, $G^i$ is simple.
\item for some $i$, $G^i$ is perfect and the group $G^i/H^{i+1}$ does not
have any proper abelian normal subgroup.
\end{itemize}}
\ms

Also the braid group on more than $4$ strings are the examples to show
that an arbitrary finite index normal subgroup of an adorable group need
not be adorable. 

In Section 4 the following result is proved about Lie groups.

\ms
\n
{\bf Theorem (4.9).} {\it Every connected real or complex Lie group is
adorable as an abstract group.}
\ms

Below we give some examples of nonadorable groups. Proofs of
nonadorability of some of these examples are easy. Proofs for the other
examples are given in the next sections.

Some examples of groups which are not adorable are nonabelian free
groups and fundamental groups of surfaces of genus greater than
$1$; for the intersection of a monotonically decreasing sequence of
characteristic subgroups of a nonabelian free group consists of the
trivial element only. The commutator subgroup of $SL_2({\Bbb Z})$ is
the nonabelian free group on 2 generators. Hence $SL_2({\Bbb Z})$ is not 
adorable. Also by Stallings' theorem if the fundamental group of a compact
$3$-manifold has finitely generated nonabelian commutator subgroup
which is not isomorphic to the Klein bottle group with infinite cyclic
abelianization then the group is not adorable. It is known that 
most of these $3$-manifolds support hyperbolic metric by Thurston.
It is easy to show that the pure braid group is not adorable as there is
a surjection of any pure braid group of more than $2$ strings onto a
nonabelian free group.

From now on whenever we give examples of nonadorable groups we will
mention its close relationship with nonpositively curved Riemannian
manifolds. This will help us to state a general conjecture (Conjecture
(0.1)).

The next result gives some important classes of examples of nonadorable
groups which are generalized free products $G_1*_HG_2$ or $HNN$-extensions
$K*_H$. We always assume $G_1\neq H\neq G_2$ and $K\neq H$.

\ms
\n
{\bf Theorem (2.3).} {\it Let $G$ be a group. 

If $G=G_1*_HG_2$ is a generalized free product and $G^1\cap H=(1)$, then
one of the following holds.

\begin{itemize}
\item $G^1$ is perfect.
\item $G^1$ is isomorphic to the infinite dihedral group.
\item $G$ is not adorable.
\end{itemize}

If $G=K*_H=\langle K, t\;\ tHt^{-1}=\phi(H)\rangle$ is an $HNN$-extension
and $G^1\cap H=(1)$, then $G$ is not adorable.

In the second case and in the last possibility of the first case for
$i\geq 1$ the rank of $G^i/G^{i+1}$ is $\geq 2$.}

\ms

In Corollary (2.7) we deduce a more general version of Theorem (2.3) and
show that if $H$ is $n$-step $G$-solvable (see Definition (2.6)) then in
the amalgamated free product case either $G$ is adorable of degree at
most $n+1$ or is not adorable and in the $HNN$-extension case it is always
nonadorable. 

We will give some more examples (Lemma (2.8) and Example (2.9)) of
a class of nonadorable generalized free products and examples of
compact Haken $3$-manifolds with nonadorable fundamental groups. 

At this point, recall that if $M$ is a connected, closed oriented
$3$-manifold and $\pi_2(M, x)\neq 0$ then by Sphere theorem
([\cite{He}, p. 40]) there is an embedded $2$-sphere in $M$ representing
a nonzero element of $\pi_2(M, x)$. Hence $M$ can be written as a
connected sum of two nonsimply connected $3$-manifolds and thus
$\pi_1(M,x)$ is a nontrivial free product. In addition if we assume that
$\pi_1(M,x)$ is not perfect and $M$ is not the connected sum of two
projective $3$-spaces then by Theorem (2.3) $\pi_1(M,x)$ is not
adorable. Thus we see that most closed $3$-manifolds with $\pi_2(M, x)\neq
0$ have nonadorable fundamental groups.

The next result is about groups with some geometric assumption. Recall
that a torsion free Bieberbach groups is the fundamental group of a
Riemannian manifold with sectional curvature equal to $0$ everywhere. 

\ms
\n
{\bf Corollary (4.3).} {\it A torsion free Bieberbach group is
nonadorable unless it is solvable.}
\ms

The following Theorem is dealing with groups under some homological
hypothesis. This theorem has an interesting application in knot theory and
possibly in $3$-manifolds in general also.

\ms
\n
{\bf Theorem (4.4).} {\it Let $G$ be a group satisfying the following
properties.
\begin{itemize}
\item $H_1(G, {\Bbb Z})$ has rank $\geq 3$. 
\item $H_2(G^j, {\Bbb Z})=0$ for $j\geq 0$.
\end{itemize}

Then $G$ is not adorable. Moreover, $G^j/G^{j+1}$ has rank $\geq 3$ for
each $j\geq 1$.}
\ms

The Proposition below is a consequence of the above Theorem.

\ms
\n
{\bf Proposition (4.7).} {\it A knot group is adorable if and
only if it has trivial Alexander polynomial.}
\ms

In fact in this case the commutator subgroup of the knot group is perfect.
All other knot groups are not adorable. On the other hand any knot
complement supports a complete nonpositively curved Riemannian metric
({\cite{L}). 

After seeing the preprint (\cite{R2}) Tim Cochran informed me that the
Proposition (4.7) was also observed by him in [\cite{C}, corollary 4.8]. 

Note that most of the torsion free examples of nonadorable groups we
mentioned above act freely and properly discontinuously (except the braid
groups case, which is still an open question) on a simply connected
complete nonpositively curved Riemannian manifold. Also we recall
that a solvable subgroup of the fundamental group of a nonpositively
curved manifold is virtually abelian (\cite{Y}). There are
generalization of these results to the case of locally $CAT(0)$ spaces
(\cite{BH}). Considering these facts we pose the following conjecture.

\ms
\n
{\bf Conjecture (0.1).} {\it Fundamental group of generic class
of complete nonpositively curved Riemannian manifolds or more generally of
generic class of locally $CAT(0)$ metric spaces are not
adorable.}
\ms

One can even ask the same question for hyperbolic groups.

Now we state the conjecture we referred before. Though in \cite{R2} this
conjecture was stated for any finitely presented torsion free groups our
primary aim was the following particular case.    

\ms
\n
{\bf Conjecture (0.2).} {\it Let $G$ be a finitely presented torsion free 
group which is isomorphic to the fundamental group of a closed
aspherical $3$-manifold such that $G^i/G^{i+1}$ is a finite group for all
$i$. Then $G$ is adorable.} 
\ms

Using Theorem (3.1) in Section 3 it is easy to show that the above
conjecture is true for aspherical Seifert fibered spaces. In fact we will
show that most Seifert fibered spaces have nonadorable fundamental groups. 

Here note that a partial converse of the above conjecture is true for 
closed $3$-manifold. Before we prove this claim note that the hypothesis
of the conjecture implies that each $G^i$ is finitely generated.

\ms
\n
{\bf Lemma (0.3).} {\it Let $G$ be the fundamental group of a closed
$3$-manifold such that for some $i$, $G^i$ is nontrivial, finitely
generated and perfect. Then for each $i$, $G^i/G^{i+1}$ is a finite
group.}

\ms
\n
\begin{proof} Since $G^i$ is a nontrivial perfect group, it is not
a surface group. Also since $G^i$ is finitely generated, by [\cite{He},
theorem 11.1] $G^i$ is of finite index in $G$. This proves the
Lemma.\end{proof}

\ms
\n
{\bf Remark (0.4).} After seeing the preprint \cite{R2} Peter A. Linnell
pointed out to me that certain finite index subgroups of $SL(n, {\Bbb Z})$
for $n\geq 3$ satisfy the hypothesis of [\cite{R2}, conjecture 0.2] but
they are not adorable. These are some noncocompact lattices in $SL(n,
{\Bbb R})$ which are residually finite $p$-groups and satisfy Kazhdan
property T. I thank Professor Linnell for the stimulating example. We
describe his example in the Appendix. Conjecture (0.2) remains open for
the fundamental groups of closed aspherical $3$-manifolds and for
cocompact discrete subgroups of Lie groups. Considering this situation we
state our main problem below.

\ms
\n
{\bf Main Problem.} {\it Find groups for which the Conjecture 0.2 is
true.}
\ms

Note that $G^i/G^{i+1}$ is finite for each $i$ if and only if $G/G^i$
is finite for each $i$. Thus, in other words the above conjecture says
that a nonadorable aspherical $3$-manifold group has an infinite
solvable quotient. Compare this observation with Proposition (4.1). 

Also note that by Theorem (2.3), if the group $G$ in Conjecture (0.2) is
not perfect and not isomorphic to ${\Bbb Z}_2*{\Bbb Z}_2$ then it is
irreducible. Thus we can assume that the group $G$ in the Conjecture is
irreducible. Recall that a group is {\it irreducible} if the group is not
isomorphic to free product of two nontrivial groups. 

There is another consequence of this conjecture. That is, if Conjecture
(0.2) is true then the virtual Betti number conjecture will be
true if a modified (half) version of it is true. We mention it below.

\ms
\n
{\bf Modified virtual Betti number conjecture.} {\it Let $M$ be a closed 
aspherical $3$-manifold such that $H_1(M, {\Bbb Z})=0$. Then there is a 
finite sheeted covering $\tilde M$ of $M$ with $H_1(\tilde M, {\Bbb Z})$ 
infinite.}
\ms

It is easy to see that the Conjecture (0.2) and the Modified virtual Betti
number conjecture together implies the virtual Betti number conjecture.

\ms
\n
{\bf Virtual Betti number conjecture.} {\it Any closed aspherical
$3$-manifold has a finite sheeted covering with infinite first homology
group.}
\ms

The virtual Betti number conjecture was raised as a question by John
Hempel in [\cite{He1}, question 1.2]. 

\section{Some elementary facts about  adorable\\ groups}

Recall that a group is called {\it perfect} if the commutator subgroup of
the group is the whole group.

\ms
\n
{\bf Lemma (1.1).} {\it Let $f:G\to H$ be a surjective homomorphism
with $G$ adorable. Then $H$ is also adorable and $doa(H)\leq
doa(G)$.}

\ms
\n
{\bf Example (1.2).} The Artin pure braid group on more than $2$ strings
is not adorable, for it has a quotient a nonabelian free group. 
In fact the full braid group on $n$-strings is not adorable for
$n\leq 4$ and adorable of degree $1$ otherwise. (see \cite{GL}).

\ms
\n
{\bf Lemma (1.3).} {\it The product $G\times H$ of two groups are
adorable if and only if both the groups $G$ and $H$ are adorable. Also if
$G\times H$ is adorable then $doa(G\times H)=\mbox{max}\ \{doa(G),
doa(H)\}$.}
\ms

On the contrary, in the case of free product of groups, almost all the
time the output is nonadorable. Hence, adorability is mainly a property
for irreducible groups. We will consider the case of free product and more
generally the generalized free product case in the next section. 

\ms
\n
{\bf Lemma (1.4).} {\it Let $G$ be an adorable group and $H$ is a normal 
subgroup of $G$. Assume that for some $i_0$, $G^{i_0}$ is simple. Then 
$H$ is also adorable and $doa(H)\leq doa(G)$.}

\ms
\n
{\bf Remark (1.5).} In the above lemma instead of assuming the strong 
hypothesis that $G^{i_0}$ is simple we can assume only that $G^{i_0}$ 
is perfect and $G^{i_0}/H^{i_0+1}$ does not have any proper normal abelian
subgroup. With this weaker hypothesis the proof follows from the fact that
the kernel of the surjective homomorphism $G^{i_0}/H^{i_0+1}\to
G^{i_0}/H^{i_0}$ is either trivial or $G^{i_0}=H^{i_0}$. In either case 
it follows that $H$ is adorable.

\ms
\n
{\bf Lemma (1.6).} {\it Let $H$ be a normal subgroup of a group $G$ such
that $G^i/H^i$ is abelian for some $i$. Then $G$ is adorable if and only
if $H$ is adorable.}

\ms
\n
{\bf Proposition (1.7).} {\it Let $H$ be a normal subgroup of a group 
$G$ such that $G/H$ is solvable. Then $H$ is adorable if and only if so is
$G$.}

\ms
\n
\begin{proof} Before we start the proof, we note down some generality. 
Suppose $G$ has a filtration of the form $G_n< G_{n-1}<\cdots <G_1<G_0=G$ 
where $G_i$ is normal in $G_{i-1}$ and $G_{i-1}/G_i$ is abelian for each
$i$. Since $G_{i-1}/G_i$ is abelian for each $i$, we have 
${G'}_{i-1}\subset G_i$. Replacing $i$ by $i+1$ we get $G_i'\subset
G_{i+1}$. Consequently, $G_0^i=G^i=\{G'\}^{i-1}\subset {G_1}^{i-1}\subset
{\{G_1}'\}^{i-2} \subset {G_2}^{i-2}\subset\cdots \subset
{G'}_{i-1}\subset G_i$. Thus we get $G^n\subset G_n$. 

Denote $G/H$ by $F$. As $F$ is solvable we have $1\subset F^k\subset\cdots 
\subset F^1\subset F^0=F$ where $F^k$ is abelian. Let $\pi:G\to G/H$ be
the quotient map. We have the following sequence of normal subgroups of
$G$. $$\cdots\subset H^n\subset H^{n-1}\cdots\subset H^1\subset H\subset
\pi^{-1}(F^k)\cdots\subset \pi^{-1}(F^0)=G.$$

Note that this sequence of normal subgroups satisfy the same properties 
as those of the filtration $G_i$ of $G$ above. Hence $G^{k+i}\subset
H^{i-1}$. Now if $G$ is adorable then for some $i$, $G^{k+i}$ is
perfect. We have $$H^{k+i}\subset G^{k+i}=G^{k+k+i+2}\subset
H^{k+i+1}.$$ But we already have $H^{k+i+1}\subset H^{k+i}$. That is
$H^{k+i}$ is perfect, hence $H$ is adorable. Conversely if $H$ is adorable
then for some $i$, $H^i$ is perfect. Note from the above inclusions that 
$H^i=G^i$ for some large $i$. Hence $G$ is also adorable.\end{proof}

\ms
\n
{\bf Theorem (1.8).} {\it A group $G$ is adorable if and only if    
there is a filtration $G_n< G_{n-1}<\cdots <G_1<G_0=G$ of
$G$ so that $G_i$ is normal in $G_{i-1}$, $G_{i-1}/G_i$ is abelian for
each $i$ and $G_n$ is a perfect group.}

\ms
\n
\begin{proof} We use Proposition (1.7) and induction on $n$ to prove
the `if' part of the Theorem. So assume that there is a filtration of $G$
as in the hypothesis. Then $G_n$ is an adorable subgroup of $G$ with
solvable quotient $G/G_n$. Proposition (1.7) proves this implication. The
`only if' part of the Theorem follows from the definition of adorable
groups.\end{proof}

\ms
\n
{\bf Corollary (1.9).} {\it Let $G$ be a torsion free infinite group and
$F$ be a finite quotient of $G$ with kernel $H$ such that $H$ is free
abelian and also central in $G$. Then $G$ is adorable.}

\ms
\n
\begin{proof} Recall that equivalence classes of extensions of $F$ by 
$H$ are in one to one correspondence with $H^2(F, H)$ which is isomorphic
to $Hom(F, ({\Bbb R}/{\Bbb Z})^n)$ where $n$ is the rank of $H$ 
([\cite{Br}, p. 95, exercise 3]). If $F$ is perfect then $Hom(F, ({\Bbb
R}/{\Bbb Z})^n)=0$ and hence the extensions $1\to H\to G\to F\to 1$
splits. But by hypothesis $G$ is torsion free. Hence $F$ is not perfect.
By a similar argument it can be shown that $F^i$ is perfect for no $i$
unless it is the trivial group. Since $F$ is finite this proves that $F$
is solvable and hence $G$ is adorable, in fact solvable.\end{proof}

We sum up the above Lemmas and Propositions in the following Theorem.

\ms
\n
{\bf Theorem (1.10).} {\it Let $H$ be a normal subgroup of an adorable
group $G$. Then $H$ is adorable if one of the following conditions is 
satisfied.
\begin{itemize}
\item $G/H$ is solvable.
\item for some $i$, $G^i/H^i$ is abelian.
\item for some $i$, $G^i$ is simple.
\item for some $i$, $G^i$ is perfect and the group $G^i/H^{i+1}$ does not
have any proper abelian normal subgroup.
\end{itemize}}

\ms
\n
{\bf Remark (1.11).} It is known that any countable group is a subgroup of
a countable simple group (see [\cite{LS}, chapter IV, theorem 3.4]). Also
we mentioned before that even finite index normal subgroup of an adorable
group need not be adorable. So the above theorem is best possible in this
regard.

\ms

In the next section we give some more examples of virtually adorable
groups which are not adorable. 

The following is an analogue of a theorem of Hirsch for poly-cyclic 
groups. The proofs of Lemmas (A) and (B) in the proof of the theorem are
easy and we leave it to the reader. 

\newpage
\n
{\bf Theorem (1.12).} {\it The following are equivalent.

\begin{itemize}
\item $G$ is a group which admits a filtration $G=G_0>G_1>\cdots >G_n$
with the property that each $G_{i+1}$ is normal in $G_i$ with quotient
$G_i/G_{i+1}$ cyclic and $G_n$ is a perfect group which satisfies the
maximal condition for subgroups.

\item $G$ is adorable and it satisfies the maximal condition for 
subgroups, i.e., for any sequence $H_1<H_2<\cdots$ of subgroups of 
$G$ there is an $i$ such that $H_i=H_{i+1}=\cdots$.

\end{itemize}}

\ms
\n
\begin{proof} The proof is on the same line as Hirsch's theorem. The main
lemma is the following.

\ms
\n
{\bf Lemma (A).} {\it Let $H_1$ and $H_2$ be two subgroup of a group $G$
and 
$H_1\subset H_2$. Let $H$ be a normal subgroup of $G$ with the 
property that $H\cap H_1=H\cap H_2$ and the subgroup generated by $H$ and
$H_1$ is equal to the subgroup generated by $H$ and $H_2$. Then
$H_1=H_2$.}

\ms  

{\bf (1) implies (2):} By Theorem (1.11) it follows that $(1)$ implies 
$G$ is adorable. Now we check the maximal condition by induction 
on $n$. As $G_n$ already satisfy maximal condition we only need to 
check that $G_{n-1}$ also satisfy maximal condition which follows from 
the following Lemma and by noting that $G_{n-1}/G_n$ is cyclic. 

\ms
\n
{\bf Lemma (B).} {\it Let $H$ be a normal subgroup of a group $G$ such
that 
both $H$ and $G/H$ satisfy the maximal condition then $G$ also satisfies
the maximal condition.}

\ms
\n
\begin{proof} Let $K_1<K_2<\cdots$ be an increasing sequence of subgroups
of $G$. Consider the two sequences of subgroups $H\cap K_1<H\cap
K_2<\cdots$ and $\{H, K_1\}<\{H, K_2\}<\cdots$. Here $\{A, B\}$ denotes 
the subgroup generated by the subgroups $A$ and $B$. As $H$ and $G/H$ both
satisfy the maximal condition there are integers $k$ and $l$ so that 
$H\cap K_k=H\cap K_{k+1}=\cdots$ and $\{H, K_l\}=\{H, K_{l+1}\}=\cdots$. 
Assume $k\geq l$. Then by Lemma (A) $K_k=K_{k+1}=\cdots$.\end{proof}

{\bf (2) implies (1):} As $G$ is adorable it has a filtration
$G=G_0>G_1>\cdots >G_n$ with $G_n$ perfect and each quotient abelian. 
Also $G_n$ satisfies maximal condition as it is a subgroup of $G$ and $G$
satisfies maximal condition. Since $G$ satisfies maximal condition 
each quotient $G_i/G_{i+1}$ is finitely generated. Now a filtration as in 
$(1)$ can easily be constructed. 

This proves the theorem.\end{proof}

\section{Generalized free products and adorable\\ groups}

We begin this section with the following result on free product of groups.

Recall that the infinite dihedral group $D_{\infty}$ is isomorphic to
${\Bbb Z}\rtimes {\Bbb Z}_2\simeq {\Bbb Z}_2*{\Bbb Z}_2$.

\ms
\n
{\bf Proposition (2.1).} {\it The free product $G$ of two nontrivial
groups,
one of which is not perfect, is either isomorphic to $D_{\infty}$ or not
adorable. Moreover, in the nonadorable case the rank of the abelian group
$G^i/G^{i+1}$ is greater or equal to $2$ for all $i\geq 1$.}

\ms
\n
\begin{proof} Let $G$ be the free product of the two nontrivial
groups $G_1$ and $G_2$ and one of $G_1$ and $G_2$ is not perfect. Then, as
the abelianization of $G=G_1*G_2$ is isomorphic to $G_1/G_1^1\oplus
G_2/G_2^1$, $G$ is also not perfect.

By Kurosh Subgroup theorem ([\cite{LS}, proposition 3.6]) any subgroup of
$G$ is isomorphic
to a free product $*_iA_i*F$, where $F$ is a free group and the groups 
$A_i$ are conjugates of subgroups of either $G_1$ or $G_2$. In particular
the commutator subgroup $G^1$ is isomorphic to $*_iA_i*F$ for some $A_i$
and $F$. Note that $[G_1,G_2]=\langle g_1g_2g_1^{-1}g_2^{-1}\ |\ g_i\in
G_i, i=1,2\rangle$ is a subgroup of $G^1$. Now assume that $G$ is not
$D_{\infty}$. Then $[G_1,G_2]$ is a nonabelian free group and clearly
$[G_1,G_2]\cap G_1=(1)=[G_1,G_2]\cap G_2$. Also $[G_1,G_2]$ is not
conjugate to any subgroup of $G_1$ or $G_2$. Hence $[G_1,G_2]$ is a
subgroup of $F$, which shows that $F$ is a nontrivial nonabelian free
group. Hence the abelianization of $G^1$ is nontrivial. By a similar
argument using Kurosh Subgroup theorem we conclude that no $G^n$ is
perfect. This proves the first assertion of the Proposition. The second
part follows from the fact that the free group $F$ has rank $\geq 2$ and
a nonabelian free group has derived series consisting of nonabelian free
groups.\end{proof}

\ms
\n
{\bf Remark (2.2).} In Proposition (2.1) we have seen that the free
product of any nontrivial group with a nonperfect group is either
$D_{\infty}$ or nonadorable group. The natural question arises here is
what happens in the amalgamated free product case of two groups along a
nontrivial group or in the case of $HNN$-extension? At first recall that
there are examples of simple groups which are amalgamated free product of
two nonabelian free groups along a (free) subgroup (see \cite{BM}). We
give another example. Let $M={\Bbb S}^3-N(k)$ be a knot complement of a
knot $k$ in the $3$-sphere. Assume that the Alexander polynomial of the
knot $k$ is nontrivial. Then by Proposition (4.7) we know that $\pi_1(M)$
is not adorable. Recall that the first homology of $M$ is generated by a
meridian of the torus boundary of $M$ and the longitude which is parallel
to the knot in ${\Bbb S}^3$ represents the zero in $H_1(M, {\Bbb Z})$. Now
glue two copies of $M$ along the boundary which sends the above longitude
of one copy to the meridian of the other and vice versa. Then the
resulting manifold $N$ has fundamental group isomorphic to the
amalgamated free product $\pi_1(M)*_{{\Bbb Z}\times {\Bbb Z}}\pi_1(M)$ and
an application of Mayer-Vietoris sequence for integral homology shows that
$N$ has trivial first homology. That is $N$ has perfect fundamental group.
Another example in this connection is the fundamental group of a torus
knot complement in ${\Bbb S}^3$. This group is of the form $G={\Bbb
Z}*_{{\Bbb Z}}{\Bbb Z}$. If the knot is of type $(p,q)$ then the two
inclusions of ${\Bbb Z}$ in ${\Bbb Z}$ in the above amalgamated free
product are defined by multiplication by $p$ and $q$ respectively. But $G$
is not adorable as it has nonabelian free commutator subgroup.
In the following theorem we consider a more general situation.
\ms

From now on, whenever we consider generalized free product
$G=G_1*_HG_2$ or $HNN$-extension $G=K*_H$, unless otherwise stated, we
always assume that $G_1\neq H\neq G_2$ and $K\neq H$.

\ms
\n
{\bf Theorem (2.3).} {\it Let $G$ be a group. 

If $G=G_1*_HG_2$ is a generalized free product and $G^1\cap H=(1)$, then
one of the following holds.

\begin{itemize}
\item $G^1$ is perfect.
\item $G^1$ is isomorphic to the infinite dihedral group $D_{\infty}$.
\item $G$ is not adorable.
\end{itemize}

If $G=K*_H=\langle K,t \mid tHt^{-1}=\phi(H)\rangle$ is an $HNN$-extension
and $G^1\cap H=(1)$, then $G$ is not adorable.

In the second case and in the last possibility of the first case for
$i\geq 1$ the rank of $G^i/G^{i+1}$ is $\geq 2$.}

\ms

Note that the assumption $G^1\cap H=(1)$ implies that $H$ is abelian.

To prove the theorem we need to recall the bipolar structure on
generalized free product and the characterization of
generalized free product by the existence of a bipolar structure on the 
group by Stallings. 

\ms
\n
{\bf Definition (2.4).} ([\cite{LS}, p. 207, definition]) A {\it
bipolar structure} on a group $G$ is a partition of $G$ into five disjoint
subsets $H, EE, EE^*, E^*E, E^*E^*$ satisfying the following axioms. (The
letters $X,Y,Z$ will stand for the letters $E$ or $E^*$ with the
convention that $(X^*)^*=X$, etc.)

\begin{itemize}
\item $H$ is a subgroup of $G$.
\item If $h\in H$ and $g\in XY$, then $hg\in XY$.
\item If $g\in XY$, then $g^{-1}\in YX$. ({\it Inverse axiom})
\item If $g\in XY$ and $f\in Y^*Z$, then $gf\in XZ$. ({\it Product axiom})
\item If $g\in G$, there is an integer $N(g)$ such that, if there exist
$g_1,\ldots, g_n\in G$ and $X_0,\ldots, X_n$ with $g_i\in X_{i-1}^*X_i$
and $g=g_1\cdots g_n$, then $n\leq N(g)$. ({\it Boundedness axiom})
\item $EE^*\neq\emptyset$. ({\it Nontriviality axiom})
\end{itemize}}
\ms

It can be shown that every amalgamated free product or $HNN$-extension has
a bipolar structure [\cite{LS}, p. 207-208]. The following theorem
of Stallings shows that the converse is also true.

\ms
\n
{\bf Theorem (2.5).} ([\cite{LS}, theorem 6.5]) {\it A group $G$ has a
bipolar structure if and only if $G$ is either a nontrivial free product
with amalgamation (possibly an ordinary free product) or an
$HNN$-extension.}

\ms
\n
{\it Proof of Theorem (2.3).} At first note that the first $5$ properties
in the above definition are hereditary, that is any subgroup $F$ of $G$
has
a partition by subsets satisfying these properties. The induced
partition of $F$ is obtained by taking intersection of $H, EE,\ldots$ with
$F$. But $EE^*\cap F$ could be empty. We replace $F$ by the commutator
subgroup $G^1$ of $G$. We would like to check the sixth property (that is,
the nontriviality axiom) for this induced partition on $G^1$. 

We consider the amalgamated free product case first. Recall that if we
write $g\in G-H$ in the form $g=c_1\cdots c_n$ where no $c_i\in H$
and each $c_i$ is in one of the factors $G_1$ or $G_2$ and successive
$c_i, c_{i+1}$ come from different factors, then $g\in EE^*$ if and only
if $c_1\in G_1$ and $c_n\in G_2$. Such a word is called {\it cyclically
reduced}. Thus $EE^*$ consists of all cyclically reduced words. Let
$g_1\in G_1-H$ and $g_2\in G_2-H$, then $g_1g_2g_1^{-1}g_2^{-1}$ is a
cyclically reduced word and is contained in $EE^*\cap G^1$. Hence the
induced partition on $G^1$ defines a bipolar structure on $G^1$ with
amalgamating subgroup $G^1\cap H=(1)$. Hence $G^1$ is a free product of
two nontrivial groups. Using Proposition (2.1) we complete the proof in
this case.

When $G$ is an $HNN$-extension we have a similar situation. We have to
check that $EE^*\cap G^1\neq\emptyset$. Recall from [\cite{LS}, p. 208]
that if
we write $g\in G-H$ in the reduced form $g=h_0t^{\epsilon_1}h_1\cdots
t^{\epsilon_n}h_n$ (where $\epsilon_i=\pm 1$ and $h_i\in K$ for each
$i$) then $g\in EE^*$ if and only if $h_0\in K-H$, or $h_o\in H$ and
$\epsilon =+1$, and $h_n\in H$ and $\epsilon_n=+1$. Now let $h_0\in K-H$
and $h_1\in H$, then
$h_0(h_1t^{-1})h_0^{-1}(h_1t^{-1})^{-1}=(h_0h_1)t^{-1}h_0^{-1}th_1^{-1}\in 
EE^*\cap G^1$. Hence the induced partition on $G^1$ gives a bipolar
structure on $G^1$. Since $G^1\cap H=(1)$ we get that $G^1$ is a free
product of a nontrivial group with the infinite cyclic group. Hence
Proposition (2.1) applies again.$\hfill{\Box}$     

We introduce below a stronger version of the notion of solvability which
depends both on the group and the group where it is embedded. 

\ms
\n
{\bf Definition (2.6).} A subgroup $H$ of a group $G$ is called {\it
$G$-solvable} (or {\it subgroup solvable}) if $G^n\cap H=(1)$ for some
$n$. If in addition $G^{n-1}\cap H\neq (1)$ then $H$ is called {\it
$n$-step $G$-solvable} (or {\it $n$-step subgroup solvable}).
\ms

Note that if $H$ is $G$-solvable then $H$ is solvable. Also if $G$ is
solvable then any subgroup of $G$ is $G$-solvable. 

Now we can state a Corollary of Theorem (2.3). The proof is easily deduced
from the proof of Theorem (2.3) and is left to the reader. 

\ms
\n
{\bf Corollary (2.7).} {\it Let $G$ be a group.

If $G=G_1*_HG_2$ is a generalized free product and $H$ is $n$-step
$G$-solvable, then one of the following holds. 

\begin{itemize}
\item $G$ is adorable of degree $n$ and not solvable.
\item $G^n\simeq D_{\infty}$.
\item $G$ is not adorable. 
\end{itemize}

If $G=K*_H=\langle K, t\;\ tHt^{-1}=\phi(H)\rangle$ is an $HNN$-extension
and $H$ is $G$-solvable, then $G$ is not adorable.

In the second case and in the last possibility of the first case for
$i\geq 1$ the rank of $G^i/G^{i+1}$ is $\geq 2$.}

\ms

The following Lemma consider some more generalized free product cases.

\ms
\n
{\bf Lemma (2.8).} {\it Let $G_1*_HG_2$ be a generalized free product with
$H$ abelian and is contained in the center of both $G_1$ and $G_2$.
Also assume that one of $G_1/H$ or $G_2/H$ is not perfect. Then
$G_1*_HG_2$ is either solvable or not adorable.}

\ms
\n
\begin{proof} Using normal form of elements of $G_1*_HG_2$ it is easy to
show that the center of $G_1*_HG_2$ is $H$. This implies that we have a
surjective homomorphism $G_1*_HG_2\to (G_1*_HG_2)/H=G_1/H*G_2/H$. By
Proposition (2.1) $G_1/H*G_2/H$ is either the infinite dihedral group or
not adorable and hence $G_1*_HG_2$ is either solvable or not adorable by
Lemma (1.1).\end{proof} 

\ms
\n
{\bf Example (2.9).} Using Lemma (2.8) we now give a large class of
examples of compact Haken $3$-manifolds with nonadorable fundamental
groups. Let $M$ and $N$ be two compact orientable Seifert fibered
$3$-manifolds with nonempty boundary and orientable base orbifold. Such
examples of $M$ and $N$ are torus knot complements in ${\Bbb S}^3$. Let
$\partial M$ and $\partial N$ be the boundary components of $M$ and $N$
respectively. Note that both $\partial M$ and $\partial N$ are tori. Let
$\gamma_1\subset\partial M$ and $\gamma_2\subset\partial N$ be simple
closed curves which are parallel to some regular fiber of $M$ and $N$
respectively. Recall that both $\gamma_1$ and $\gamma_2$ represent central
elements of $\pi_1(M)$ and $\pi_1(N)$ respectively. Now choose an annulus
neighborhood $A_1$ of $\gamma_1$ in $\partial M$ and $A_2$ of $\gamma_2$
in $\partial N$ and glue $M$ and $N$ identifying $A_1$ with $A_2$ by a
diffeomorphism which sends $\gamma_1$ to $\gamma_2$. Let $P$ be the
resulting manifold. Then $P$ is a compact Haken $3$-manifold with tori
boundary and by Seifert-van Kampen theorem $\pi_1(P)$ satisfies the
hypothesis of Lemma (2.8) and hence either solvable or not adorable. Here
note that the manifold $P$ itself is Seifert fibered. In the next section
we will show that in fact an infinite group which is the fundamental group
of a compact Seifert fibered $3$-manifold is nonadorable except for some
few cases.

\ms

\section{Adorability and 3-manifolds}

Seifert fibered spaces are a fundamental and very important class of
$3$-manifolds. Conjecturally (due to Thurston) any $3$-manifold is build
from Seifert fibered spaces and hyperbolic $3$-manifolds. Results of
Jaco-Shalen, Johannson and Thurston say that this is in fact true for any
Haken $3$-manifold.  

\ms
\n
{\bf Theorem (3.1).} {\it Let $M^3$ be a compact Seifert fibered
$3$-manifold. Then one of the following four cases occur.

\begin{itemize}
\item $(\pi_1(M))^i$ is finite for some $i\leq 2$.
\item $\pi_1(M)$ is solvable.
\item $\pi_1(M)$ is not adorable and $(\pi_1(M))^i/(\pi_1(M))^{i+1}$ has
rank greater than $1$ for all $i$ greater than some $i_0$.
\item $\pi_1(M)$ is perfect.
\end{itemize}}

\ms
\n
\begin{proof} At first we recall some well known group theoretic
informations about the fundamental group of Seifert fibered spaces. 
If $B$ is the base orbifold of $M$ then there is
a surjective homomorphism $\pi_1(M)\to \pi_1^{orb}(B)$, where
$\pi_1^{orb}(B)$ is the orbifold fundamental group of $B$. Recall that 
$\pi_1^{orb}(B)$ is a Fuchsian group. Also recall that the
above surjective homomorphism is part of the following exact sequence.
$$1\to \langle t\rangle \to \pi_1(M)\to \pi_1^{orb}(B)\to 1.$$ Here
$\langle t\rangle$ is the cyclic normal subgroup of $\pi_1(M)$ generated
by a regular fiber of the Seifert fibration of $M$. Also if $\pi_1(M)$ is
infinite then $\langle t\rangle$ is an infinite cyclic subgroup of
$\pi_1(M)$.

Some examples of Seifert fibered $3$-manifolds with finite fundamental
group are lens spaces and the Poincare sphere. So, from now on we assume
$\pi_1(M)$ is infinite. Then there is an exact sequence. $$1\to {\Bbb
Z}\to \pi_1(M)\to \pi^{orb}_1(B)\to 1.$$  There are now two cases to
consider. 

\n
{\bf Case 1.} $\pi^{orb}_1(B)$ is finite. By [\cite{FJ}, lemma 2.5]
$\pi_1(M)$ has a finite normal subgroup $G$ with quotient isomorphic
either to $\Bbb Z$ or to $D_{\infty}$. Since $D_{\infty}$ is solvable
$(\pi_1(M))^i$ is finite for some $i\leq 2$.

\n
{\bf Case 2.} $\pi^{orb}_1(B)$ is infinite and not a perfect group.
Then by [\cite{S}, theorem 1.5] there is a torsion free normal subgroup
$H$ of $\pi^{orb}_1(B)$ so that $\pi^{orb}_1(B)/H$ is a finite solvable
group. Hence by Proposition (1.7) $\pi^{orb}_1(B)$ is adorable if and only
if so is $H$. Since $H$ is of finite index in $\pi^{orb}_1(B)$ by a result
of Hoare, Karrass and Solitar [\cite{LS}, chapter III, proposition 7.4]
$H$ is again a
Fuchsian group. But a torsion free Fuchsian group is the fundamental group
of a compact surface (evident from the presentation of such groups). Hence
$H$ is either $\Bbb Z$ or ${\Bbb Z}\times {\Bbb Z}$ or ${\Bbb Z}\rtimes
{\Bbb Z}$ or nonadorable. Thus by Proposition (1.7) $\pi^{orb}_1(B)$ is
either solvable or nonadorable. If $\pi^{orb}_1(B)$ is solvable then from
the above exact sequence it follows that $\pi_1(M)$ is also solvable. On
the other hand Lemma (1.1) shows $\pi_1(M)$ is nonadorable whenever
$\pi^{orb}_1(B)$ is. 

Next, consider the case when $\pi^{orb}_1(B)$ is a perfect group.
Let $x_1,x_2,\ldots , x_n$ be the cone points on $B$ with indices $p_1,
p_2,\ldots , p_n$ greater than or equal to $2$. By [\cite{S}, theorem 1.5]
$\pi^{orb}_1(B)$ is perfect if and only if $B={\Bbb S}^2$ and the indices
$p_1, p_2,\ldots , p_n$ are pairwise coprime. It is well known that in
this situation $M$ is an integral homology $3$-sphere and hence $\pi_1(M)$
is also perfect. This proves the theorem.\end{proof}

Notice that the proof of the above theorem is not very illuminating in the 
sense that it does not show the cases when the groups are nonadorable or
solvable. Below we show that in fact in most cases the fundamental group
of a compact Seifert fibered space is nonadorable. For simplicity of
presentation we consider Seifert fibered spaces whose base orbifold
$B$ is orientable and has only cone singularities. Note that the proof of
the Theorem deals with both orientable and nonorientable $B$ and for any
kind of singularities. At first let us consider the case when $M$ has
nonempty boundary. Since $B$ also has nonempty boundary, $\pi^{orb}_1(B)$
is a free product of cyclic groups (\cite{He}) and hence by Proposition
(2.1) $\pi^{orb}_1(B)$ is either the infinite dihedral group or is
nonadorable if it is a nontrivial free product. Hence either $\pi_1(M)$ is
solvable (when $\pi^{orb}_1(B)$ is dihedral or cyclic) or (by Lemma (1.1))
$\pi_1(M)$ is not adorable. 

If $M$ is closed then we have the same situation as above except that
$\pi^{orb}_1(B)$ has the following form. $$\pi^{orb}_1(B)=\langle
a_1,\ldots , a_g, b_1,\ldots , b_g, x_1,\ldots , x_n\ |\
x_1^{j_1}=\cdots =x_n^{j_n}=1 ;$$$$
\Pi_{j=1}^g [a_j, b_j]x_1\cdots x_n=1\rangle$$ where $x_1,\ldots
x_n$ represents loops around cone points of $B$. We will consider the case  
$g=0$ at the end of the proof. If $g\geq 1$ then adding the extra
relations $a_1=1$ we get that $\pi^{orb}_1(B)$ has the following
homomorphic image $$\langle a_2,\ldots , a_g, b_1,\ldots , b_g, x_1,\ldots
, x_n\ |\ x_1^{j_1}=\cdots =x_n^{j_n}=1 ; $$$$\Pi_{j=2}^g [a_j,
b_j]x_1\cdots
x_n=1\rangle.$$ If there is no cone point on $B$ and $g=1$ then $M$ is an
${\Bbb S}^1$-bundle over the torus and hence has solvable fundamental
group. Otherwise the last group is a free product of the infinite cyclic
group (generated by $b_1$) and another group and hence not adorable by
Proposition (2.1). Thus $\pi^{orb}_1(B)$ is also not adorable by Lemma
(1.1). Consequently so is $\pi_1(M)$. 

Now we consider the case when $g=0$. There are further two cases to
consider.

\n
{\bf Case A.} $\pi^{orb}_1(B)$ is finite. This case occurs when $B$ has
at most 3 cone points and if exactly 3 cone points with indices
$n_1,n_2,n_3$ then $\frac{1}{n_1}+\frac{1}{n_2}+\frac{1}{n_3}>1$ (see
[\cite{He}, theorem 12.2]). We have already discussed this case in {\bf
Case 1} in the proof of the theorem. 

\n
{\bf Case B.} $\pi^{orb}_1(B)$ is infinite. In this case there
are the following two possibilities (see [\cite{He}, theorem 12.2]). $(a)$
$B$ has more than $3$ cone points. $(b)$ $B$ has $3$ cone points with
indices $j_1, j_2, j_3$ so that $\frac{1}{j_1} + \frac{1}{j_2} +
\frac{1}{j_3}\leq 1$.

For $(a)$ we need the following easily verified remark.

\ms
\n
{\bf Remark (3.2).} If $B$ is a sphere with $3$ cone points then
$|\pi_1^{orb}(B)|\geq 3$.
\ms

Now recall that in $(a)$ $\pi_1^{orb}(B)$ has the following presentation.
$$\langle x_1,\ldots , x_n\ |\ x_1^{j_1}=\cdots =x_n^{j_n}=1 ;
x_1\cdots x_n=1\rangle$$ where $n\geq 4$. Now assume $n\geq 6$ and add
the relation $x_1x_2x_3=1$ in the above presentation. Then
$\pi_1^{orb}(B)$ surjects onto the free product of $$\langle
x_1,x_2,x_3\ |\ x_1^{j_1}=x_2^{j_2}=x_3^{j_3}=1 ; x_1x_2 x_3=1\rangle$$
and 
$$\langle x_4,\ldots , x_n\ |\ x_4^{j_4}=\cdots =x_n^{j_n}=1 ; x_4\cdots
x_n=1\rangle.$$ By
Proposition (2.1) and Remark (3.2) it follows that $\pi_1^{orb}(B)$ is
either perfect or not adorable and hence so is $\pi_1(M)$. In the case
$n=5$ if there is a pair of indices $j_k$ and $j_l$ so that $(j_k,
j_l)\geq 3$ then it is easy to show that $\pi_1(M)$ is nonadorable. We
leave the remaining cases to the reader. 

In $(b)$ when $\frac{1}{j_1} + \frac{1}{j_2} + \frac{1}{j_3}=1$ then
$\pi_1^{orb}(B)$ is a discrete group of isometries of the Euclidean plane.
Recall that a torsion free discrete group of isometries of the Euclidean
plane is isomorphic to $\Bbb Z$ or ${\Bbb Z}\times {\Bbb Z}$ or ${\Bbb
Z}\rtimes {\Bbb Z}$ and hence by the result of Sah we mentioned above
$\pi_1^{orb}(B)$ is either perfect or solvable. On the other hand if
$\frac{1}{j_1} + \frac{1}{j_2} + \frac{1}{j_3}<1$ then $\pi_1^{orb}(B)$ is
a discrete groups of isometries of the hyperbolic plane. Since a group of
isometries of the hyperbolic plane does not contain a free abelian group
on more than one generator, it follows by the result of Sah that in this
case $\pi_1^{orb}(B)$ is either perfect or a finite solvable extension of
$\Bbb Z$ or nonadorable. Hence $\pi_1(M)$ is either solvable or perfect or
nonadorable.

\ms
\n
{\bf Remark (3.3).} Recall that a Fuchsian group is a discrete subgroup
of $PSL(2,{\Bbb R})$ and it is either a free product of cyclic groups or
is isomorphic to a group of the form $\pi^{orb}_1(B)$. In the free product
case except for the infinite dihedral group all other free products are
nonadorable. In the remaining cases we have already seen in the proof of
Theorem (3.1) that a Fuchsian group is either finite or perfect or
solvable or nonadorable and in most cases it is nonadorable. It is not
known to me if a similar situation occur for discrete subgroups of
$PSL(2,{\Bbb C})$. Such informations will be very useful to get some hold
on the virtual Betti number conjecture for hyperbolic
$3$-manifolds.
\ms

\section{(Non)adorability under homological or geometric hypothesis}

In Section 2 under some group theoretic hypothesis we showed
when a generalized free product or an $HNN$-extension produces
a nonadorable group. 

This section deals with some homological or geometric (or
topological) hypothesis on a group which ensures that the group is
nonadorable.

\ms
\n
{\bf Proposition (4.1).} {\it Let $M^3$ be a compact $3$-manifold with the 
property that there is an exact sequence of groups $1\to H\to \pi_1(M)\to
F\to 1$ such that $H$ is finitely generated nonabelian but not the
fundamental group of the Klein bottle and $F$ is an
infinite solvable group. Then $\pi_1(M)$ is not adorable.}

\ms
\n
\begin{proof} By theorem 11.1 in \cite{He} it follows that $H$ is the
fundamental group of a compact surface. Also as $H$ is nonabelian and not
the Klein bottle group, it is not adorable. The Proposition now follows
from Proposition (1.7).\end{proof}

\ms
\n
{\bf Proposition (4.2).} {\it Let $G$ be a torsion free group and $H$ a
free nonabelian (or abelian) normal subgroup of $G$ with quotient $F$ a
nontrivial finite (or finite perfect) group. Then $G$ is not
adorable.}

\ms
\n
\begin{proof} If $H$ is nonabelian then by Stallings' theorem $G$
itself is free and hence not adorable. So assume $H$ is free abelian. 
Since in this case $F$ is a perfect group, the restriction of the quotient
map $G\to F$ to $G^i$ is again surjective for each $i$ with $H\cap G^i$ as
kernel. And since $G$ is infinite and torsion free, $H\cap G^i$ is
nontrivial free abelian for all $i$. This shows that each $G^i$ is
again a Bieberbach group. Note that if $H^1(G^i, {\Bbb Z})=0$ then 
$G^i$ is centerless and it is known that centerless Bieberbach groups are
meta-abelian with nontrivial abelian holonomy group and hence solvable
(\cite{HS}). But since each $G^i$ surjects 
onto a nontrivial perfect group it cannot be solvable. Hence 
$H^1(G^i, {\Bbb Z})\neq 0$ for each $i$. This proves the
Proposition.\end{proof}

The conclusion of the above Proposition remains valid if we assume that
$F$ is nonsolvable adorable.

By Bieberbach theorem (\cite{Ch}) we have the following Corollary.

\ms
\n
{\bf Corollary (4.3).} {\it The fundamental group of a closed flat
Riemannian manifold is nonadorable unless it is solvable.}
\ms

So far we have given examples of nonadorable groups which are
fundamental groups of known class of manifolds or of manifolds with some 
strong Riemannian structure. The following Theorem gives a general class
of examples of nonadorable groups under some homological conditions.

\ms
\n
{\bf Theorem (4.4).} {\it Let $G$ be a group satisfying the following
properties.
\begin{itemize}
\item $H_1(G, {\Bbb Z})$ has rank $\geq 3$. 
\item $H_2(G^j, {\Bbb Z})=0$ for $j\geq 0$.
\end{itemize}

Then $G$ is not adorable. Moreover, $G^j/G^{j+1}$ has rank $\geq 3$ for
each $j\geq 1$.}

\ms
\n
\begin{proof} Consider the short exact sequence. $$1\to G^1\to G\to
G/G^1\to 1.$$

We use the Hochschild-Serre spectral sequence ([\cite{Br}, p. 171]) of
the above exact sequence. The $E^2$-term of the spectral sequence is  
$E_{pq}^2=H_p(G/G^1, H_q(G^1, {\Bbb Z}))$. Here ${\Bbb Z}$ is
considered as a trivial $G$-module. This spectral sequence gives rise to
the following five term exact sequence. $$H_2(G, {\Bbb Z})\to
E^2_{20}\to E^2_{01}\to H_1(G, {\Bbb Z})\to E^2_{10}\to 0.$$ 

Using $(2)$ we get $$0\to H_2(G/G^1, H_0(G^1, {\Bbb Z}))\to H_0(G/G^1,
H_1(G^1, {\Bbb Z}))\to H_1(G, {\Bbb Z})$$$$\to H_1(G/G^1, H_0(G^1,
{\Bbb Z}))\to 0.$$ 

As $\Bbb Z$ is a trivial $G$-module we get $$0\to H_2(G/G^1, {\Bbb Z})\to
H_0(G/G^1, H_1(G^1, {\Bbb Z}))\to H_1(G, {\Bbb Z})\to $$$$H_1(G/G^1, {\Bbb 
Z})\to 0.$$

Note that the homomorphism between the last two nonzero terms in the
above exact sequence is an isomorphism. Also the second nonzero term from
left is isomorphic to the co-invariant $H_1(G^1, {\Bbb Z})_{G/G^1}$ and
hence we have the following $$H_2(G/G^1, {\Bbb Z})\simeq H_1(G^1, {\Bbb
Z})_{G/G^1}.$$ 

Since $G/G^1$ has rank $\geq 3$ we get that $H_2(G/G^1, {\Bbb Z})$ has
rank greater or equal to $3$. This follows from the following lemma. 

\ms
\n
{\bf Lemma (4.5).} {\it Let $A$ be an abelian group. Then the rank
of $H_2(A, {\Bbb Z})$ is $rk A(rk A -1)/2$ if $rk A$ is finite
and infinite otherwise.}

\ms
\n
\begin{proof} If $A$ is finitely generated then from the formula
$H_2(A, {\Bbb Z})\simeq \bigwedge^2 A$ it follows that rank of $H_2(A,
{\Bbb Z})$ is $rk A(rk A -1)/2$. In the case $A$ is countable and
infinitely generated then there are finitely generated subgroups $A_n$ of
$A$ such that $A$ is the direct limit of $A_n$. Now as homology of group
commutes with direct limit the proof follows using the previous
case. Similar argument applies when $A$ is uncountable.\end{proof}

To complete the proof of the theorem note that there is a surjective
homomorphism $H_1(G^1, {\Bbb Z})\to H_1(G^1, {\Bbb Z})_{G/G^1}$. Thus we
have proved that $H_1(G^1, {\Bbb Z})$ also has rank $\geq 3$. Finally
replacing $G$ by $G^n$ and $G^1$ by $G^{n+1}$ and using induction on $n$
the proof is completed.\end{proof}

There are two important consequences of Theorem (4.4). At first we recall
some definition from \cite{St}. 

Let $R$ be a nontrivial commutative ring with unity. The class $E(R)$
consists of groups $G$ for which the trivial $G$-module $R$ has a
$RG$-projective resolution $$\cdots\to P_2\to P_1\to P_0\to R\to 0$$ 
such that the map {\bf$1$}$_R\otimes \partial_2:R\otimes_{RG}P_2\to
R\otimes_{RG}P_1$ is injective. Note that if a group belongs to  $E(R)$
then $H_2(G, R)=0$. Also this condition is sufficient to belong to
$E(R)$ for groups of cohomological dimension less or equal to $2$. By
definition $G$ lies in $E$ if it belongs to $E(R)$ for all $R$. A 
characterization of $E$-groups is that a group $G$ is an $E$-group if and
only if $G$ belongs to $E({\Bbb Z})$ and $G/G^1$ is torsion free
([\cite{St}, lemma 2.3]).

\ms
\n
{\bf Corollary (4.6).} {\it Let $G$ be an $E$-group and rank of $H_1(G,
{\Bbb Z})$ is $\geq 2$. Then $G$ is not adorable.}

\ms
\n
\begin{proof} By [\cite{St}, theorem A] it follows that $G$ satisfies
the second condition of Theorem (4.4). Hence we get  that $H_1(G^2, {\Bbb
Z})$ has rank $\geq 1$ and hence in particular $G^2$ is not perfect. On
the other hand an $E$-groups has derived length $0,1,2$ or infinity
(remark after [\cite{St}, theorem A]). Thus $G$ is not 
adorable.\end{proof}

In the following Proposition we give an application of Theorem (4.4) 
for knot groups.

\ms
\n
{\bf Proposition (4.7).} {\it Let $H=\pi_1({\Bbb S}^3-k)$, where $k$ is a 
nontrivial knot in the $3$-sphere with nontrivial Alexander polynomial.
Then $H$ is not adorable. Moreover if rank of $H^1/H^2$ is
greater or equal to $3$ then the same is true for
$H^j/H^{j+1}$ for all $j\geq 2$.}
\ms

In fact a stronger version of the Proposition follows, namely by \cite{St}  
the successive quotients of the derived series of $G$ are
torsion free. Thus we get that the successive quotients of the derived
series are nontrivial and torsion free. 

\ms
\n
{\it Proof of Proposition (4.7).} At first recall that the second 
condition of Theorem (4.4) follows from [\cite{St}, theorem A]. On the
other hand the commutator subgroup of a knot group is perfect if and only
if the knot has trivial Alexander polynomial. So assume that $H^1$ is not
perfect. If $H^1$ is finitely generated then in fact it is nonabelian
free and hence $H$ is not adorable. If rank of $H^1/H^2$ is $\geq 3$ then
the proof follows from the above Theorem. So assume that rank of $H^1/H^2$
is $\leq 2$. 

Recall that the rank of the abelian group $H^1/H^2$ is equal to the degree
of the Alexander polynomial of the knot (see [\cite{Cr}, theorem 1.1]).
Thus if rank of $H^1/H^2$ is $1$ then the Alexander polynomial has degree
$1$ which is impossible as the Alexander polynomial of a knot always has
even degree. Next if rank of $H^1/H^2$ is $2$ then $H$ is not adorable by
Corollary (2.6) and noting that knot groups are $E$-groups.$\hfill{\Box}$

\ms
\n
{\bf Definition (4.8).} A Lie group is called {\it adorable} if it
is adorable as an abstract group.

\ms
\n
{\bf Theorem (4.9).} {\it Every connected (real or complex) Lie group is
adorable.}

\ms
\n
\begin{proof} Let $G$ be a Lie group and consider its derived 
series. $$\cdots\subset G^n\subset G^{n-1}\cdots \subset G^1\subset
G^0=G.$$ 

Note that each $G^i$ is a normal subgroup of $G$. Define
$G_i=\overline{G^i}$. Then we have a sequence of normal subgroups 
$$\cdots\subset G_n\subset G_{n-1}\cdots \subset G_1\subset G_0=G$$ 
so that $G_i$ is a closed Lie subgroup of $G$ and $G_i/G_{i+1}$ is
abelian for each $i$. Suppose for some $i$, dim $G_i=0$,  
i.e., $G_i$ is a closed discrete normal subgroup of $G$. We claim 
$G_i$ is abelian. For, fix $g_i\in G_i$ and consider the continuous map
$G\to G_i$ given by $g\mapsto gg_ig^{-1}$. As $G$ is connected 
and $G_i$ is discrete image of this map is the singleton $\{g_i\}$. That
is $g_i$ commutes with all $g\in G$ and hence $G_i$ is abelian. 

As $G^i\subset G_i$, $G^i$ is also abelian. Thus $G$ is solvable and hence 
adorable.

Next assume no $G_i$ is discrete. Then as $G$ is finite dimensional and
$G_i$'s are Lie subgroups of $G$ there is an $i_0$ so that $G_j=G_{j+1}$
for all $j\geq i_0$ and dim $G_{i_0}\geq 1$. We need the following Lemma
to complete the proof of the Theorem.

\ms
\n
{\bf Lemma (4.10).} {\it Let $G$ be a (real or complex) Lie group
such that $\overline {G^1}=G$. Then $G^2=G^1$, that is $G^1$ is a perfect
group.}

\ms
\n
\begin{proof} The proof of the lemma follows from [\cite{Ho}, theorem
XII.3.1 and theorem XVI.2.1].\end{proof}

We have $G^{i_0}\subset G_{i_0}$ and hence $$G_{i_0}=G_{i_0+1}
=\overline{G^{i_0+1}}\subset \overline{G^1_{i_0}}\subset
\overline{G_{i_0}}=G_{i_0}.$$ This implies $\overline{G^1_{i_0}}=G_{i_0}$. 
Now from the above Lemma we get $G_{i_0}$ is adorable. Thus $G_{i_0}$ is a
normal adorable subgroup of $G_{i_0-1}$ with quotient $G_{i_0-1}/G_{i_0}$ 
abelian and hence by Proposition (1.7) $G_{i_0-1}$ is also adorable. By
induction it follows that $G$ is adorable.\end{proof}

\section{Appendix}

In this section we describe the counter example given by Peter A. Linnell
to [\cite{R2}, conjecture 0.2]. 

\ms
\n
{\bf Example (5.1).} (P.A. Linnell) Let $n\geq 3$ and $p$ be an odd prime.
Let $K$ be the kernel of the homomorphism $SL(n, {\Bbb Z})\to SL(n, {\Bbb
Z}/p{\Bbb Z})$ which is induced by the homomorphism ${\Bbb Z}\to {\Bbb
Z}/p{\Bbb Z}$. When $p=2$, let $K$ be the kernel of $SL(n, {\Bbb
Z})\to SL(n, {\Bbb Z}/4{\Bbb Z})$. Now we have the following three facts
about $K$.

\begin{itemize}
\item $K$ is a residually finite $p$-group. Hence we get that $K^{i+1}$ is
a proper subgroup of $K^i$ for each $i$.

\item $K$ satisfies Kazhdan property $T$. Hence $K^i/K^{i+1}$ is a finite
group for each $i$.

\item $K$ is finitely presented and torsion free.
 
\end{itemize}

Thus $K$ is not adorable. But by the second and the third fact above, $K$
satisfies the hypothesis of [\cite{R2}, conjecture (0.2)].

\ms

A notable fact is that $K$ is a noncocompact discrete subgroup of $SL(n,
{\Bbb R})$. It will be very interesting to prove Conjecture (0.2) for
cocompact discrete subgroup of Lie groups.

\section{Problems}

In this section we state some problems for a further study on adorable
groups. We also give the motivations behind each problem and mention known
results related to the problem.

\ms
\n
{\bf Problem (6.1).} {\it Study the Main Problem for some particular class
of groups, for example for cocompact discrete subgroups of Lie groups or
for groups which are fundamental groups of closed nonpositively curved
Riemannian manifolds.}

\ms

Problem (6.1) is related to the particular case of the virtual Betti
number conjecture for hyperbolic $3$-manifolds. We have already seen that
a discrete subgroup of $PSL(2, {\Bbb R})$ is either finite or solvable or
perfect or nonadorable. In fact it is possible to describe when each of
these possibilities occur. A similar result about discrete subgroup of
$PSL(2, {\Bbb C})$ will be very important. A more precise problem is the
following.

\ms
\n
{\bf Problem (6.2).} {\it Given a positive integer $n$ does there exist
a discrete (torsion free) subgroup of $PSL(2, {\Bbb C})$ which is adorable
of degree $n$?}

\ms
\n
{\bf Problem (6.3).} {\it Find all $3$-manifolds with adorable fundamental
group.}

\ms

Some examples of such $3$-manifolds are integral homology $3$-spheres and 
knot complement of knots with trivial Alexander polynomial. In Theorem
(3.1) we have seen that most Seifert fibered spaces have nonadorable
fundamental group and also we have shown when the fundamental group is
adorable.

\ms
\n
{\bf Problem (6.4).} {\it Prove that most groups are not adorable.}

\ms

A possible approach to study Problem (6.4) is by the same method which was
used to show that most groups are hyperbolic.

A small and first step towards Conjecture (0.2) is the following.

\ms
\n
{\bf Problem (6.5).} {\it Show that Conjecture (0.2) is true for the
fundamental groups of compact Haken $3$-manifolds.}

\ms

We have already mentioned that it is true for Seifert fibered spaces. Note
that if the fundamental group of a compact Haken $3$-manifold satisfies
the hypothesis of Conjecture (0.2) then the manifold has to be closed.

\vspace{.5cm}
\begin{flushleft}
{\it Address:}\\
School of Mathematics\\
Tata Institute\\
Homi Bhabha Road\\
Mumbai 400 005, India.\end{flushleft}

\begin{flushleft}
{\it email address:}-roushon@math.tifr.res.in\\
{\it homepage URL:}-http://www.math.tifr.res.in/\~\
roushon/paper.html\end{flushleft}

\newpage
\bibliographystyle{plain}
\ifx\undefined\bysame
\newcommand{\bysame}{\leavevmode\hbox to3em{\hrulefill}\,}
\fi

\enddocument